\documentclass[a4paper,reqno,11pt]{amsart}
\usepackage{amsfonts}

\theoremstyle{plain}
\begingroup
\newtheorem{theorem}{Theorem}[section]
\newtheorem{lemma}[theorem]{Lemma}
\newtheorem{proposition}[theorem]{Proposition}

\endgroup

\theoremstyle{definition}
\begingroup

\newtheorem{remark}[theorem]{Remark}

\endgroup

\theoremstyle{remark}

\mathchardef\emptyset="001F

\numberwithin{equation}{section}

\newcommand{\R}{\mathbb R}

\newcommand{\leb}{{\mathcal L}}
\newcommand{\mthree}{{\mathbb M}^{3{\times}3}}
\newcommand{\mskw}{{\mathbb M}^{3{\times}3}_{skw}}
\newcommand{\dist}{{\rm dist}}
\newcommand{\sym}{{\rm sym}\,}
\newcommand{\skw}{{\rm skew}\,}
\renewcommand{\div}{{\rm div}}
\newcommand{\wto}{\rightharpoonup}
\newcommand{\eps}{\varepsilon}

\newcommand{\E}{{\mathcal E}}

\title[Convergence of equilibria of thin elastic beams]
{Convergence of equilibria of three-dimensional thin elastic beams}
\author{M.G.\ Mora}
\author{S.\ M\"uller}
\address[M.G.~Mora]{SISSA, Via Beirut 2-4, 34014 Trieste, Italy}
\address[S.~M\"uller]{Max-Planck Institute for Mathematics in the Sciences, 
Inselstrasse 22, 04103 Leipzig, Germany}
\email[Maria Giovanna Mora]{mora@sissa.it}
\email[Stefan M\"uller]{sm@mis.mpg.de}

\begin{document}
\begin{abstract}
A convergence result is proved for the equilibrium configurations
of a three-dimensional thin elastic beam, as the diameter $h$ of the cross-section goes to zero.
More precisely, we show that stationary points of the nonlinear elastic functional $E^h$,
whose energies (per unit cross-section) are bounded by $C h^2$, converge to stationary
points of the $\Gamma$-limit of $E^h/h^2$.
This corresponds to a nonlinear one-dimensional model for inextensible rods,
describing bending and torsion effects.
The proof is based on the rigidity estimate for low-energy
deformations by Friesecke, James, and M\"uller \cite{FJM02}
and on a compensated compactness argument in a singular
geometry. In addition, possible concentration effects of the strain are controlled
by a careful truncation argument.
\end{abstract}

\maketitle

{\small

\bigskip

}

\bigskip
\bigskip

\section{Introduction and main result}

In this paper we extend our previous work with M.G.\ Schultz on the convergence of equilibria
of planar thin elastic beams (see \cite{MMS06}) to the case of three-dimensional thin beams.

To set the stage let $h>0$ and let $S$ be a bounded open connected
subset of $\R^2$ with Lipschitz boundary. We consider a thin beam whose reference
configuration is given by the open set $\Omega_h= (0,L){\times}hS$.
Given any deformation $v\in W^{1,2}(\Omega_h;\R^3)$, we define the elastic energy
(per unit cross-section) associated to $v$ as
$$
E^h (v) := \frac{1}{h^2} \int_{\Omega_h} W( \nabla v) \, dz.
$$
The stored-energy density function $W:\mthree\to [0,+\infty]$ is assumed to
satisfy the following conditions:
\begin{itemize}
\item[(h1)] frame indifference:
$W(RF) = W(F)$ for every $R\in SO(3)$ and $F\in\mthree$;
\smallskip
\item[(h2)] $W=0$ on $SO(3)$;
\smallskip
\item[(h3)] $W(F) \geq c\,  \dist^2(F,SO(3))$, $c>0$, for every $F\in\mthree$;
\smallskip
\item[(h4)] $W$ is of class $C^2$ in a neighbourhood of $SO(3)$. 
\end{itemize}
Here $SO(3)$ denotes the group of proper rotations. 
The frame indifference implies that there exists a function $\tilde{W}$
defined on symmetric matrices such that $W(\nabla v) = \tilde{W}((\nabla v)^T \nabla v)$;
i.e., the elastic energy depends only on the pull-back metric of $v$.

To discuss the limiting behaviour of $E^h$, as $h \to 0$, it is convenient to rescale
to a fixed domain $\Omega= (0,L){\times}S$ by the change
of variables
$$
z = (x_1, h x_2, h x_3) \quad\hbox{and}\quad  y(x) = v(z(x)).
$$ 
With the notation
$$
\nabla_h y = (\partial_1 y| \tfrac1h \partial_2 y | \tfrac1h \partial_3 y) 
$$
we can write the elastic energy as
$$
E^h(v) = I^h(y) := \int_\Omega W(\nabla_h y) \, dx.
$$

Without loss of generality we can assume that $\leb^2(S)=1$ and that the 
segment $(0,L){\times}\{0\}{\times}\{0\}$ is a line of centroids for the beam; i.e.,
\begin{equation}\label{section}
\int_S x_2\, dx_2dx_3=\int_S x_3\, dx_2dx_3=
\int_S x_2x_3\, dx_2dx_3 =0.
\end{equation}
Under the previous assumptions it is possible to identify a complete hierarchy of
limiting rod theories, depending on the scaling of $I^h$,
by means of $\Gamma$-convergence.
More precisely, for every $\beta\geq0$ we have
\begin{equation}\label{Gamma}
\frac{1}{h^\beta}I^h \stackrel{\Gamma}{\longrightarrow} I_\beta,
\end{equation}
where, according to $\beta$, the functional $I_\beta$ describes a different elastic
model for rods.
The $\Gamma$-convergence for $\beta=0$ was 
proved by Acerbi, Buttazzo, and Percivale in \cite{ABP91}, leading to
a {\em nonlinear string model\/}. The scaling $\beta=2$, which corresponds to a 
{\em nonlinear rod model\/}, has been studied in \cite{Mor-Mue} and
independently by Pantz in \cite{Pa02}. 
The result for $\beta=4$ has been proved in \cite{Mor-Mue2},
while for the other scalings $\Gamma$-convergence can be easily derived
from \cite{Mor-Mue} and \cite{Mor-Mue2}.

The $\Gamma$-convergence results \eqref{Gamma} guarantee that
if $(y^{(h)})$ is a compact sequence of minimizers of $I^h$ (with respect
to some boundary conditions or body forces) such that $I^h(y^{(h)})\leq Ch^\beta$,
then, up to subsequences, $(y^{(h)})$ converges to a minimizer of $I_\beta$
(for a comprehensive introduction to $\Gamma$-convergence we refer to \cite{DM}).
 
In this paper we deal with the problem of the convergence of {\em equilibria} in
the scaling $\beta=2$. In this case
the natural class of admissible functions for the limit problem turns out to be 
\begin{multline}
\mathcal{A} := \left\{ (y,d_2,d_3) \in W^{2,2}((0,L); \R^3){\times} W^{1,2}((0,L);\R^3)
{\times} W^{1,2}((0,L);\R^3): 
\right.
\\
\left.
R:=( y' | d_2 | d_3)\in SO(3) \ \hbox{a.e.\ in } (0,L)
\right\}.
\end{multline}
On this class the $\Gamma$-limit functional is given by
\begin{equation}\label{I_2}
I_2(y,d_2,d_3) := \frac{1}{2} \int_0^L  Q_1 (R^T R') \, dx_1,
\end{equation}
where $R:=(y' | d_2 | d_3)$. The density $Q_1$ is a quadratic
form on the space $\mskw$ of skew-symmetric matrices, defined as
\begin{equation}\label{minpb}
Q_1(A):=\min_{\alpha\in W^{1,2}(S;\R^3)}
\int_S Q_3\Big( x_2Ae_2+x_3Ae_3
\Big| \, \partial_2 \alpha \, \Big| \,\partial_3\alpha
\Big)\, dx_2dx_3
\end{equation}
for every $A\in \mskw$, where $Q_3$ is the quadratic form 
$Q_3(F):=\leb F{\,:\,}F$ and $\leb$ is the linear map on $\mthree$ given by
$\leb:=D^2 W(Id)$.

In this limit model the function $y$ represents the deformation of the mid-fiber of the rod,
which has to be isometric because of the constraint $|y'|=1$. The two Cosserat vectors
$d_2$ and $d_3$ determine the rotation undergone by the cross-section of the rod at each
point of the mid-fiber. We remark also that, as $R$ belongs to $SO(3)$ a.e., 
the matrix $R^T R'$ is skew-symmetric. Moreover,
the entries $(R^T R')_{1j}$ for $j=2,3$ are related to the
curvature of the deformed mid-fiber, while $(R^T R')_{23}$ is related to the torsion of the mid-fiber
and to the twist of the cross-section, after the deformation. Finally, 
the solutions to \eqref{minpb} with $A$ replaced by $R^T(x_1)R'(x_1)$
describe the warping of the cross-section with respect to the normal plane
(see \cite{Mor-Mue}).

If, in addition, $W$ is isotropic and $S$ is a disc, then 
the quadratic form $Q_1$ can be explicitly computed and reduces to
$$
Q_1(A):=\frac{1}{2\pi}\frac{\mu(3\lambda+2\mu)}{\lambda+\mu}(A_{12}^2+A_{13}^2)
+\frac{\mu}{2\pi}A_{23}^2,
$$
where $\lambda$ and $\mu$ are the Lam\'e coefficients of the rod (see 
\cite[Remark~3.5]{Mor-Mue}).

We assume the beam to be subject to a body force of density $h^2g$, with $g\in L^2((0,L);\R^3)$;
thus, we consider the functionals
\begin{equation}
J^h(y) = \int_\Omega ( W(\nabla_h y)  - h^2 g(x_1){\,\cdot\,} y ) \, dx.
\end{equation}
The corresponding $\Gamma$-limit at scale $h^2$ is then given by
\begin{equation}
J_2(y,d_2,d_3) = I_2(y,d_2,d_3)- \int_0^L g {\,\cdot\,} y \, dx_1
\end{equation}
if $(y,d_2,d_3) \in \mathcal{A}$, while
$J_2$ takes the value $+ \infty$ if $(y,d_2,d_3) \notin \mathcal{A}$
(here we took the liberty to identify maps on $\Omega$ which
are independent of $x_2,x_3$ with maps on $(0,L)$).
It is convenient to fix one end of the rod by requiring, e.g., $y(0) = 0$ and $d_k(0)=e_k$ for $k=2,3$,
where $\{e_1,e_2,e_3\}$ denotes the canonical basis in $\R^3$.

We are now in a position to state the main theorem of the paper.

\begin{theorem}\label{thm:1.1}
Assume that (h1)--(h4) are satisfied and that $W$ is differentiable
with globally Lipschitz derivative $DW$. Let $g\in L^2((0,L);\R^3)$.
Let $(y^{(h)})$ be a sequence of stationary points of $J^h$, subject to the boundary condition
$y^{(h)}(0,x_2,x_3) = (0,hx_2,hx_3)$ at $x_1=0$ and to
natural boundary conditions on the remaining
boundaries. Assume further that there exists a constant $C>0$ such that
\begin{equation} \label{eq:energy}
\int_\Omega W(\nabla_h y^{(h)})\,dx \leq C h^2
\end{equation}
for every $h$.
Then, up to subsequences,
\begin{eqnarray}
& \label{th1}
y^{(h)} \to \bar y \quad \mbox{in } W^{1,2}(\Omega; \R^3),
\\
& \label{th1.1}
\frac1h \partial_k y^{(h)} \to \bar d_k  \quad \mbox{in } L^2(\Omega; \R^3), \quad
k=2,3,
\end{eqnarray}
where $(\bar y, \bar d_2,\bar d_3)\in{\mathcal A}$ is a stationary point of 
$$
J_2(y,d_2,d_3)=\frac12 \int_0^L  Q_1 (R^T R') \, dx_1-  \int_0^L g {\,\cdot\,} y \, dx_1
$$
with respect to the boundary conditions $y(0)=0$, $d_k(0)=e_k$ for $k=2,3$, and natural
boundary conditions at $x_1=L$.
\end{theorem}

\begin{remark}
An easy application of the Poincar\'e inequality shows that the estimate
(\ref{eq:energy}) holds automatically for minimizers.
\end{remark}

\begin{remark}
In \cite{Mie} Mielke used a centre manifold approach to
compare solutions in a thin strip to a $1d$ problem. 
His approach gives a comparison already for finite $h$,
but it requires that the nonlinear strain $(\nabla_h y)^T \nabla_h y$
is close to the identity in $C^{0,\alpha}$ (and applied forces $g$
cannot be easily included). 
\end{remark}

In the case of planar thin beams the Euler-Lagrange equations
corresponding to the limit functional $J_2$ can be expressed in terms of a single ODE
in the variable $\theta$, describing the angle of the tangent vector to the
deformed mid-fiber with respect to a fixed direction. One of the major differences
in the case of three-dimensional thin beams is that the limiting Euler-Lagrange 
equations involve both a linear system of PDEs in the cross-section and
a system of ODEs in terms of the bending moments of the rod (see Section~2).
This requires an extra work in all the derivation argument.

However, the main ingredients of the proof of Theorem~\ref{thm:1.1} remain basically the same
as in the planar case discussed in \cite{MMS06}.
First the quantitative rigidity estimate in \cite{FJM02}
is used to define suitable strain-like and stress-like variables
$G^{(h)}$ and $E^{(h)}$, which are almost curl-free and divergence-free (see Steps~2 and 3).
Then we can argue in the spirit of the theory of compensated compactness,
developed by Murat and Tartar \cite{Mur,T1,T2}, to obtain strong compactness of the stress $E^{(h)}$.
This then allows us to pass to the limit in the Euler-Lagrange equations (see Step~7). 

To rule out possible concentration effects of the strain
a careful truncation argument for gradients in thin domains
is employed (see Lemma~\ref{truncate}). 
We emphasize that in the planar case
this result can be proved using a simple extension argument by successive reflection,
while in the $3d$ case an appropriate choice of the extension operator is needed.

\section{Preliminary results}

The aim of this section is to derive the Euler-Lagrange equations for
the functional $J_2$ introduced in the previous section.

We begin by collecting some properties of the minimum problem \eqref{minpb} 
defining the limit density $Q_1$. Using Korn's inequality and the direct method
of the calculus of variations it is easy to see that problem \eqref{minpb} has a solution.
Moreover, there exists a unique minimizer belonging to the class
$$
{\mathcal B}:=\Big\{ \alpha\in W^{1,2}(S): \ \int_S\alpha\, dx_2dx_3= 
\int_S\partial_2\alpha\, dx_2dx_3=\int_S\partial_3\alpha\,dx_2dx_3=0\Big\}
$$
(see \cite[Remark~3.4]{Mor-Mue}).
The Euler-Lagrange equations for problem \eqref{minpb} are computed in the next lemma.

\begin{lemma}\label{lm:EulerQ}
Let $A\in\mskw$ and let $F_A:W^{1,2}(S;\R^3)\to [0,+\infty)$ be the functional defined by
\begin{equation}
F_A(\alpha):= \int_S Q_3\Big(x_2Ae_2+x_3Ae_3
\Big| \, \partial_2 \alpha \, \Big| \,\partial_3\alpha
\Big)\, dx_2dx_3
\end{equation}
for every $\alpha\in W^{1,2}(S;\R^3)$. 
Then a function $\alpha\in {\mathcal B}$ is the minimizer of $F_A$ if and only if
the function $E:S\to\mthree$ given by
$$
E:=\leb\Big(x_2Ae_2+x_3Ae_3
\Big| \, \partial_2 \alpha \, \Big| \,\partial_3\alpha
\Big)
$$
satisfies (in a weak sense) the boundary value problem
\begin{equation}\label{eqlm0}
\begin{cases}
\div_{\!x_2,x_3} (Ee_2\, |\, Ee_3)=0 & \hbox{in } S,
\\
(Ee_2\,|\,Ee_3)\,\nu_{\partial S}=0 & \hbox{on }\partial S,
\end{cases}
\end{equation}
where $\nu_{\partial S}$ is the outer unit normal to $\partial S$.
Moreover, the minimizer depends linearly on the entries of $A$.
\end{lemma}

\begin{proof}
As $F_A$ is a convex functional, a function $\alpha\in{\mathcal B}$ minimizes $F_A$ if and only if it satisfies 
$$
\int_S \leb\Big(x_2Ae_2+x_3Ae_3
\Big| \, \partial_2 \alpha \, \Big| \,\partial_3\alpha
\Big){\,:\,} \Big(0\, |\, \partial_2 \beta\, |\, \partial_3 \beta \Big)\, dx_2dx_3=0
$$
for every $\beta\in W^{1,2}(S;\R^3)$, which is equivalent to \eqref{eqlm0}.
The linear dependence of $\alpha$ on the entries of $A$ follows directly from the equation
\eqref{eqlm0}.
\end{proof}

\begin{remark}
Let $(y, d_2,d_3)\in{\mathcal A}$, let $R:=(y'\,|\,d_2\,|\,d_3)$, and let $A:=R^T R'$.  For every
$x_1\in(0,L)$ let $\alpha(x_1,\cdot)\in {\mathcal B}$ be the minimizer of
\eqref{minpb} with $A$ replaced by $A(x_1)$. Since $\alpha(x_1,\cdot)$ depends linearly on $A(x_1)$
by Lemma~\ref{lm:EulerQ} and $A\in L^2((0,L);\mskw)$, we conclude that
$\alpha\in L^2(\Omega;\R^3)$ and $\partial_k \alpha\in L^2(\Omega;\R^3)$ for $k=2,3$.
\end{remark}

The next lemma is concerned with the derivation of the Euler-Lagrange equations
for the functional $J_2$. The stationary condition for
a triple $(y, d_2,d_3)\in{\mathcal A}$ satisfying the
boundary conditions will be expressed in terms of the
bending moments $\tilde E$ and $\hat E$ defined below.
Let $E:\Omega\to\mthree$ be the {\em stress} corresponding to the
deformation $(y, d_2,d_3)$, defined by
\begin{equation}\label{1dstress}
E(x):= \leb \Big(x_2A(x_1)e_2+x_3A(x_1)e_3
\Big| \, \partial_2 \alpha(x) \, \Big| \,\partial_3\alpha(x)
\Big),
\end{equation}
where $A:=R^T R'$, $R:=(y'\,|\,d_2\,|\,d_3)$, and 
$\alpha\in L^2((0,L);\R^3)$ is such that $\alpha(x_1,\cdot)\in{\mathcal B}$ solves
\eqref{minpb}, with $A$ replaced by $A(x_1)$, for a.e.\ $x_1\in(0,L)$. 
We call the {\em bending moments} associated with the deformation $(y, d_2,d_3)$ 
the functions
$\tilde E:(0,L)\to\mthree$ and $\hat E:(0,L)\to\mthree$ given by
$$
\tilde E(x_1):=\int_S x_2E(x)\, dx_2dx_3, \qquad \hat E(x_1):=\int_S x_3E(x)\, dx_2dx_3
$$
for every $x_1\in(0,L)$.

\begin{lemma}\label{ELJ2}
Let $(y, d_2,d_3)\in{\mathcal A}$ be such that $y(0)=0$ and $d_k(0)=e_k$ for $k=2,3$.
Then $(y, d_2,d_3)$ is a stationary point of $J_2$ with respect to the boundary conditions $y(0)=0$ and $d_k(0)=e_k$ for $k=2,3$ (and natural boundary conditions at $x_1=L$) if and only if
the following system of equations is satisfied:
\begin{equation}\label{eq:EL2}
\left\{
\begin{array}{l}
\tilde E_{11}'  =  A_{13}(\hat E_{21} -\tilde E_{31}) - A_{23} \hat E_{11}
-R^T\tilde g{\, \cdot\,}e_2,
\smallskip
\\
\hat E_{11}'  =  -A_{12}(\hat E_{21} -\tilde E_{31}) + A_{23} \tilde E_{11}
-R^T\tilde g{\, \cdot\,}e_3,
\smallskip
\\
\hat E_{21}' - \tilde E_{31}' =  A_{12}\hat E_{11}-A_{13}\tilde E_{11},
\smallskip
\\
\tilde E_{11}(L)=\hat E_{11}(L)=0,\quad \hat E_{21}(L)- \tilde E_{31}(L)=0,
\end{array}
\right.
\end{equation}
where 
$$
\tilde g(x_1):=\int_L^{x_1}g(t)\, dt
$$
for every $x_1\in(0,L)$.
\end{lemma}

\begin{remark}
If $W$ is isotropic, that is,
$$
W(F)=W(FR) \quad \text{for every } F\in\mthree,\, R\in SO(3),
$$ 
then the linear operator $\leb$ associated with the second derivatives of $W$ 
at the identity reduces to
$$
\leb F=2\mu\,\sym F +\lambda({\rm tr} F) Id,
$$
where $\lambda$ and $\mu$ are the Lam\'e coefficients of the rod.

If we assume in addition that the cross-section $S$ is a disc, then
the minimizer $\alpha\in{\mathcal B}$ of \eqref{minpb} can be explicitly computed
and, in terms of the entries of the matrix $A=R^TR'$,
it is given by
$$
\begin{array}{c}
\alpha = -\tfrac14\tfrac{\lambda}{\lambda+\mu}(x_2^2A_{12} -x_3^2 A_{12}
+2x_2x_3A_{13})e_2
\smallskip
\\
\hphantom{-x_2^2A_{12}} {-}\tfrac14\tfrac{\lambda}{\lambda+\mu}(-x_2^2A_{13} +x_3^2 A_{13}+2x_2x_3A_{12})e_3,
\end{array}
$$
(see \cite[Remark~3.5]{Mor-Mue}).
In this case the stress is equal to
$$
E=\left(
\begin{array}{ccc}
\frac{\mu(3\lambda+2\mu)}{\lambda+\mu}(x_2A_{12}+x_3A_{13}) &
\frac12 x_3A_{23} & -\frac12 x_2A_{23} \\
\frac12 x_3A_{23} & 0 & 0 \\
-\frac12 x_2A_{23} & 0 & 0
\end{array}
\right),
$$
while the bending moments are
\begin{eqnarray*}
\tilde E & = & \tfrac{1}{4\pi}
\tfrac{\mu(3\lambda+2\mu)}{\lambda+\mu}A_{12}e_1\otimes e_1-
\tfrac{1}{8\pi} A_{23} (e_1\otimes e_3+e_3\otimes e_1),
\\
\hat E & = & \tfrac{1}{4\pi}
\tfrac{\mu(3\lambda+2\mu)}{\lambda+\mu}A_{12}e_1\otimes e_1+
\tfrac{1}{8\pi} A_{23} (e_1\otimes e_2+e_2\otimes e_1).
\end{eqnarray*}
\end{remark}

\begin{proof}[Proof of Lemma~\ref{ELJ2}]
Let $R:=(y'\,|\,d_2\,|\,d_3)$ and let $A:=R^T R'$.
It is convenient to consider $J_2$ as a functional defined on the class
$$
{\mathcal R}:=\{P\in W^{1,2}((0,L);\mthree):\ P\in SO(3) \text{ a.e.\ in } (0,L), \ P(0)=Id \},
$$ 
whose tangent space at $R$ is given by all functions of the form $RB$
with $B\in W^{1,2}((0,L);\mskw)$ and $B(0)=0$. 

Let then $B\in W^{1,2}((0,L);\mskw)$ with $B(0)=0$. 
In order to compute the G\^ateaux differential of $J_2$ at $R$ in the tangent direction given by $RB$, 
we consider a smooth curve $\gamma:[0,1]\to {\mathcal R}$
such that $\gamma(0)=R$ and $\dot\gamma(0)=RB$ (where the dot denotes derivative with
respect to the variable $\eps\in[0,1]$). Then we have
\begin{equation}\label{minmal}
J_2(\gamma(\eps))=
\frac12 \int_0^L Q_1(\gamma(\eps)^T\gamma(\eps)')\,dx_1
+ \int_0^L \tilde g {\,\cdot\,} \gamma(\eps)e_1\, dx_1,
\end{equation}
where the prime denotes derivative with respect to $x_1\in[0,L]$.
Now, let $\beta^\eps \in L^2(\Omega;\R^3)$ be such that $\beta^\eps (x_1,\cdot)\in{\mathcal B}$
is the solution to the problem \eqref{minpb} with $A$ replaced by $\gamma(\eps)^T\gamma(\eps)'$
for a.e.\ $x_1\in(0,L)$.
Then
\begin{eqnarray*}
\lefteqn{\frac12 \int_0^L Q_1(\gamma(\eps)^T\gamma(\eps)')\,dx_1}
\\
& = &
\frac12 \int_\Omega Q_3\Big(x_2\,\gamma(\eps)^T\gamma(\eps)'e_2+x_3\,\gamma(\eps)^T\gamma(\eps)'e_3
\Big| \, \partial_2 \beta^\eps \, \Big| \,\partial_3\beta^\eps
\Big)\, dx.
\end{eqnarray*}
Differentiating equation \eqref{minmal} at $\eps=0$ and taking into account the previous
formula, we obtain
\begin{eqnarray*}
\lefteqn{dJ_2(R)[RB]}
\\
& = &
\int_\Omega E{\,:\,} \Big(x_2(AB-BA+B')e_2 + x_3(AB-BA+B')e_3 \, \Big| \, \partial_2 \beta \, \Big| \, 
\partial_3 \beta  \Big)\, dx
\\
& & + \int_0^L R^T\tilde g {\,\cdot\,} Be_1\, dx_1,
\end{eqnarray*}
where $E$ is the stress defined in \eqref{1dstress}
and $\beta \in L^2(\Omega;\R^3)$ is such that $\beta(x_1,\cdot)\in{\mathcal B}$ is the solution
to the problem \eqref{minpb} with $A$ replaced by $B^TR'+R^TB'$ for a.e.\ $x_1\in(0,L)$. 
Here we used the fact that  by Lemma~\ref{lm:EulerQ}
the function $\beta^\eps$ depends linearly on the entries of $\gamma(\eps)^T\gamma(\eps)'$.

By \eqref{eqlm0} the vectorfield $Ee_2,Ee_3$ is divergence free in the variables $x_2,x_3$,
hence
$$
\int_S ( Ee_2{\,\cdot\,} \partial_2 \beta + Ee_3{\,\cdot\,}
\partial_3 \beta  )\, dx_2dx_3= 0.
$$
Thus the differential of $J_2$ reduces to
\begin{eqnarray*}
dJ_2(R)[RB] & = &
\int_\Omega Ee_1{\,\cdot \,} (x_2(AB-BA+B')e_2 + x_3(AB-BA+B')e_3)\, dx
\\
& & + \int_0^L R^T\tilde g {\,\cdot\,} Be_1\, dx_1.
\end{eqnarray*}
Integration with respect to $x_2,x_3$ in the first term on the right-hand side yields
\begin{multline}\label{eq:EL}
dJ_2(R)[RB]= \int_0^L ( \tilde Ee_1{\,\cdot\,}B'e_2 +  \hat Ee_1{\,\cdot\,}B'e_3)\, dx_1
\\
+\int_0^L ( \tilde Ee_1{\,\cdot\,} (AB-BA)e_2 + \hat Ee_1{\,\cdot\,} (AB-BA)e_3 )\, dx_1
+\int_0^L R^T\tilde g{\,\cdot\,}Be_1\, dx_1.
\end{multline}
As $A$ is skew-symmetric, we have that for any $F\in\mthree$ and for $k=2,3$
$$
Fe_1{\,\cdot\,}(AB-BA)e_k=
-AFe_1{\,\cdot\,}Be_k -\sum_{j\neq k}A_{jk}Fe_1{\,\cdot\,}Be_j.
$$
Using \eqref{eq:EL} and the previous formula, it is easy to
check that the condition 
$$
dJ_2(R)[RB]=0 \quad \text{for every } B\in W^{1,2}((0,L);\mskw) \text{ with } B(0)=0
$$
is equivalent to the following three equations:
\begin{equation}\label{system}
\begin{array}{c}
\displaystyle\int_0^L \!\! \big( \phi'\,\tilde E_{11}
+\phi\, A_{13}(\hat E_{21}-\tilde E_{31})-\phi\, A_{23}\hat E_{11}-\phi\, 
R^T\tilde g{\,\cdot\,}e_2\big) \, dx_1=0,
\\
\displaystyle\int_0^L \!\! \big( \phi'\,\hat E_{11}
- \phi\, A_{12}(\hat E_{21}-\tilde E_{31})+\phi\, A_{23}\tilde E_{11}-\phi\, 
R^T\tilde g{\,\cdot\,}e_3\big) \, dx_1=0,
\\
\displaystyle\int_0^L \!\! \big( \phi'(\hat E_{21}-\tilde E_{31})
+ \phi\, A_{12}\hat E_{11} -\phi\, A_{13}\tilde E_{11}\big) \, dx_1=0
\end{array}
\end{equation}
for every $\phi\in W^{1,2}(0,L)$ with $\phi(0)=0$.
By integration by parts the previous equations are equivalent to system \eqref{eq:EL2}.
\end{proof}

\section{Proof of Theorem~\ref{thm:1.1}}

Let $(y^{(h)})$ be a sequence of stationary points of $J^h$; i.e., suppose that
the following condition is satisfied:
\begin{equation}\label{EL1}
\int_\Omega \Big( DW(\nabla_h y^{(h)}){\,:\,}\nabla_h\psi -h^2g\cdot\psi\Big)\, dx =0
\end{equation}
for every $\psi\in W^{1,2}(\Omega;\R^3)$ such that $\psi(0,x_2,x_3)=0$ for $(x_2,x_3)\in S$. 
Assume that (\ref{eq:energy}) holds true.

The proof is split into several steps. 
\bigskip

\noindent
{\em Step 1. Decomposition of the deformation gradients in rotation and strain.}
\smallskip

\noindent
By Proposition~\ref{rotation} there exists a sequence $(R^{(h)})\subset C^\infty((0,L);\mthree)$
such that $R^{(h)}(x_1)\in SO(3)$ for every $x_1\in (0,L)$ and
\begin{eqnarray}
&
\|\nabla_h y^{(h)}-R^{(h)}\|_{L^2}\le Ch, \label{rig1}
\\
&
\|(R^{(h)})'\|_{L^2}+h \|(R^{(h)})''\|_{L^2}\le C, \label{rig2}
\\
&
|R^{(h)}(0)-Id|\le C\sqrt{h}. \label{rig3}
\end{eqnarray}
By (\ref{rig2}), up to subsequences, $R^{(h)}$ converge to some $R$
weakly in $W^{1,2}((0,L);\mthree)$, hence uniformly in $L^\infty((0,L);\mthree)$. 
Thus $R(x_1)\in SO(3)$ for every $x_1\in(0,L)$. 
{}From inequality (\ref{rig1}) it follows that
$$
\nabla_h y^{(h)}\to R \quad \mbox{strongly in } L^2(\Omega;\mthree).
$$
In particular, we have that $\partial_k y^{(h)}\to 0$ for $k=2,3$ and thus 
\begin{equation}\label{gradient}
\nabla y^{(h)}\to Re_1\otimes e_1  \quad \mbox{strongly in } L^2(\Omega;\mthree).
\end{equation}
As $|y^{(h)}(0,x_2,x_3)|\le Ch\to 0$, we deduce from the Poincar\'e inequality that
$y^{(h)}$ converge to some $\bar{y}$ strongly in $W^{1,2}(\Omega;\R^3)$ 
and that $\bar{y}$ satisfies 
$$
\partial_1\bar{y}=Re_1, \quad 
\partial_2\bar{y}=\partial_3\bar{y}=0  \quad
\hbox{a.e.\ in }\Omega.
$$ 
Therefore, setting $\bar d_k:= Re_k$ for $k=2,3$, we have that $(\bar y,\bar d_2,\bar d_3)\in{\mathcal A}$, and the convergence properties (\ref{th1}) and (\ref{th1.1}) are proved.
Moreover, the boundary conditions at $x_1=0$ follow from \eqref{rig3} and the uniform convergence of $R^{(h)}$.

Let $G^{(h)}:\Omega\to\mthree$ be the function
$$
G^{(h)}:=\tfrac1h ((R^{(h)})^T\nabla_h y^{(h)}-Id).
$$
As the functions $G^{(h)}$ are bounded in $L^2(\Omega;\mthree)$ by (\ref{rig1}),
we can assume, up to extracting a subsequence, that
\begin{equation}\label{Gdeb}
G^{(h)}\wto G \quad \mbox{weakly in } L^2(\Omega;\mthree)
\end{equation}
for some $G\in L^2(\Omega;\mthree)$. 
Moreover, from the definition of $G^{(h)}$ it follows immediately that the deformation gradients can be decomposed as
\begin{equation}\label{decomp}
\nabla_h y^{(h)}=R^{(h)}(Id+hG^{(h)}).
\end{equation}
\medskip

\noindent
{\em Step 2. Consequence of compatibility for the strain.}
\smallskip

\noindent
The decomposition \eqref{decomp} suggests that, roughly speaking, the
strains $G^{(h)}$ have the structure of a scaled gradient, up to the factor $(R^{(h)})^T$.
This implies that the limit strain $G$ has to satisfy some compatibility constraints.
In order to deduce these conditions we introduce a sequence of auxiliary deformations
$z^{(h)}:\Omega\to\R^3$ defined by
\begin{equation}\label{def:z^h}
z^{(h)}(x):=\tfrac1h  y^{(h)}(x) -\tfrac1h \int_0^{x_1} R^{(h)}(s)e_1\, ds
- x_2R^{(h)}(x_1)e_2 - x_3R^{(h)}(x_1)e_3.
\end{equation}
Using (\ref{decomp}) we obtain
\begin{eqnarray}
\nabla_h z^{(h)} & = & \tfrac1h (\nabla_h y^{(h)} -R^{(h)}) 
- x_2(R^{(h)})'e_2 \otimes e_1 - x_3(R^{(h)})'e_3
\otimes e_1
\nonumber
\\
& = & R^{(h)}(G^{(h)}- x_2A^{(h)}e_2\otimes e_1 - x_3A^{(h)}e_3\otimes e_1),
\label{form1}
\end{eqnarray}
where $A^{(h)}:=(R^{(h)})^T(R^{(h)})'$.
Since $R^{(h)}\wto R$ in $W^{1,2}((0,L);\mthree)$, we have
\begin{equation}\label{weakA}
A^{(h)}\wto A:=R^TR' \quad \hbox{weakly in } L^2((0,L);\mthree).
\end{equation}
Using these two facts, together with \eqref{Gdeb}, we conclude that
\begin{equation}\label{nablaz2}
\nabla_h z^{(h)} \wto R (G - x_2Ae_2\otimes e_1 - x_3Ae_3\otimes e_1)
\quad \mbox{weakly in } L^2(\Omega;\mthree).
\end{equation}
As $|z^{(h)}(0,x_2,x_3)|\leq C\sqrt{h}$ by \eqref{rig3}, we deduce from the Poincar\'e inequality
that $z^{(h)}$ converge to some $z$ weakly in $W^{1,2}(\Omega;\R^3)$. 
Moreover, the limit function $z$ satisfies
\begin{equation}\label{form2}
R^T\partial_1 z= Ge_1- x_2Ae_2 - x_3Ae_3 , 
\quad \partial_2 z=\partial_3 z= 0 \quad \hbox{a.e.\ in } \Omega.
\end{equation}
In particular, $z$ does not depend on $x_2,x_3$ and, thus, by the first equality in (\ref{form2}) 
$Ge_1$ is an affine function of $x_2,x_3$.
If we denote by $\bar G$ the zeroth moment of $G$ defined by
$$
\bar G(x_1):=\int_S G(x)\, dx_2dx_3, \qquad x_1\in(0,L),
$$
then it follows immediately from \eqref{section} and \eqref{form2} that
\begin{equation}\label{Gbar}
\bar G e_1= R^T z'.
\end{equation}

To identify the second and third column of $G^{(h)}$ it is convenient
to define $\alpha^{(h)}:\Omega\to\R^3$ as 
$$
\alpha^{(h)}:=\tfrac1h (R^{(h)})^T z^{(h)} - \int_S \tfrac1h (R^{(h)})^T z^{(h)}\, dx_2dx_3.
$$
{}From (\ref{nablaz2}) and the uniform convergence of $R^{(h)}$ it follows that
\begin{equation}\label{alpha}
\partial_k \alpha^{(h)}\wto Ge_k \quad \hbox{weakly in } L^2(\Omega;\R^3)
\end{equation}
for $k=2,3$. By the Poincar\'e inequality on the cross-section $S$, there exists a constant $C>0$ such that for a.e.\ $x_1\in(0,L)$
$$
\|\alpha^{(h)}(x_1,\cdot)\|_{L^2(S)}^2\leq C\|\partial_2\alpha^{(h)}(x_1,\cdot)\|_{L^2(S)}^2
+C\|\partial_3\alpha^{(h)}(x_1,\cdot)\|_{L^2(S)}^2.
$$
Integrating with respect to $x_1$, we deduce by \eqref{alpha} that $\alpha^{(h)}\wto\alpha$ weakly in
$L^2(\Omega;\R^3)$, where $\alpha$ satisfies $\alpha\in L^2(\Omega; \R^3)$,
$\partial_k \alpha\in L^2(\Omega; \R^3)$ for $k=2,3$, and
$Ge_k=\partial_k \alpha$ for $k=2,3$.
In particular, the function 
$$
\beta(x):=\alpha(x)-x_2\int_S \partial_2\alpha\, dx_2dx_3
-x_3\int_S \partial_3\alpha\, dx_2dx_3
$$
satisfies $\beta\in L^2(\Omega; \R^3)$,
$\partial_k \beta\in L^2(\Omega; \R^3)$ for $k=2,3$,
$\beta(x_1,\cdot)\in{\mathcal B}$ for a.e.\ $x_1\in(0,L)$, and
\begin{equation}\label{form2bis}
Ge_k-\bar Ge_k=\partial_k \beta \quad \hbox{for } k=2,3.
\end{equation}
\medskip

\noindent
{\em Step 3. Consequences of the Euler-Lagrange equations.}
\smallskip

\noindent
Let $E^{(h)}:\Omega\to\mthree$ be the scaled stress defined by
\begin{equation}\label{Ehdef}
E^{(h)}:=\tfrac1h DW(Id+h G^{(h)}).
\end{equation}
Since $DW$ is Lipschitz continuous and the $G^{(h)}$ are bounded
in $L^2(\Omega;\mthree)$, the functions $E^{(h)}$ are also bounded in $L^2(\Omega;\mthree)$.
In fact, by Proposition~\ref{Ehweak} we have that
\begin{equation}\label{lebG}
E^{(h)}\wto E:=\leb G \quad \mbox{weakly in }L^2(\Omega;\mthree).
\end{equation}
We note in particular that $E$ is symmetric, as $\leb F=(\leb F)^T$ for every $F\in\mthree$.
Note also that $\leb F=\leb(\sym F)$ for every $F\in\mthree$.

By the decomposition (\ref{decomp}) and by frame indifference we obtain that
$$
DW(\nabla_h y^{(h)})= R^{(h)} DW(Id+hG^{(h)}) = h R^{(h)}E^{(h)}.
$$
Using this identity we can write the Euler-Lagrange equations (\ref{EL1})
in terms of the stresses $E^{(h)}$. More precisely, we have
\begin{equation}\label{EL2}
\int_\Omega ( R^{(h)}E^{(h)}{\,:\,}\nabla_h\psi -hg\cdot\psi )\, dx =0
\end{equation}
for every $\psi\in W^{1,2}(\Omega;\R^3)$ with $\psi=0$ on $\{x_1=0\}$. 
Multiplying (\ref{EL2}) by $h$ and passing to the limit as $h\to 0$, we get
\begin{equation}\label{cons1}
\int_\Omega ( REe_2\cdot\partial_2\psi + REe_3\cdot\partial_3\psi )\, dx =0.
\end{equation}
As $R$ is a pointwise rotation depending only on $x_1$, the previous equation
yields
\begin{equation}\label{divfree}
\begin{cases}
\div_{\!x_2,x_3} (Ee_2\, |\, Ee_3)=0 & \hbox{in } S,
\\
(Ee_2\,|\,Ee_3)\,\nu_{\partial S}=0 & \hbox{on }\partial S
\end{cases}
\end{equation}
for a.e.\ $x_1\in(0,L)$.
This implies in particular that for a.e.\ $x_1\in(0,L)$
\begin{equation}\label{Eek}
\int_S Ee_k\, dx_2dx_3=0 \quad \hbox{for } k=2,3.
\end{equation}
\medskip

\noindent
{\em Step 4. Symmetry properties of $E^{(h)}$.}
\smallskip

\noindent
{}From the frame indifference of $W$ it follows that the matrix $DW(F)F^T$ is symmetric.
Applying this with $F=Id+hG^{(h)}$, we obtain that
\begin{equation}\label{symm}
E^{(h)}-(E^{(h)})^T= -h ( E^{(h)}(G^{(h)})^T-G^{(h)}(E^{(h)})^T ).
\end{equation}
As $E^{(h)}$ and $G^{(h)}$ are bounded in $L^2(\Omega;\mthree)$,
we deduce in particular the estimate
\begin{equation}\label{cons3}
\|E^{(h)}- (E^{(h)})^T\|_{L^1}\le Ch. 
\end{equation}
\medskip

\noindent
{\em Step 5. Moments of the Euler-Lagrange equations.}
\smallskip

\noindent
Let us introduce the zeroth and first moments of the stress $E^{(h)}$,
defined by
$$
\bar E^{(h)}(x_1):=\int_S E^{(h)}(x)\,dx_2dx_3,
$$
$$
\tilde E^{(h)}(x_1):=\int_S x_2E^{(h)}(x)\,dx_2dx_3,
 \qquad
\hat E^{(h)}(x_1):=\int_S x_3E^{(h)}(x)\,dx_2dx_3
$$
for every $x_1\in (0,L)$. We shall derive the Euler-Lagrange equations satisfied 
by the moments.

Let $\varphi\in C^\infty([0,L];\R^3)$ be such that $\varphi(0)=0$. Using $\varphi$ as
test function in the Euler-Lagrange equation (\ref{EL2}), we obtain
$$
\int_\Omega ( R^{(h)}E^{(h)}e_1{\,\cdot\,}\varphi' -hg{\,\cdot\,}\varphi)\, dx =0.
$$
Integrating first with respect to the variables of $S$ and taking into account that
$R^{(h)}$, $\varphi$, and $g$ depend only on the $x_1$ variable, we can rewrite the previous
equality as
$$
\int_0^L ( R^{(h)}\bar E^{(h)}e_1{\,\cdot\,}\varphi' -hg{\,\cdot\,}\varphi )\, dx_1 =0.
$$
Since this equation holds for every $\varphi\in C^\infty([0,L];\R^3)$ with $\varphi(0)=0$,
we deduce that
\begin{equation}\label{mom0}
\bar E^{(h)}e_1= -h(R^{(h)})^T \tilde g \quad \hbox{a.e.\ in } (0,L),
\end{equation}
where $\tilde g$ is the primitive of $g$ defined in (\ref{eq:EL}). In particular, passing
to the limit, we obtain
\begin{equation}\label{mom1}
\bar Ee_1=0 \quad \hbox{a.e.\ in } (0,L).
\end{equation}
Together with \eqref{Eek}, this implies that $\bar E=0$ a.e.\ in $(0,L)$.
As $E=\leb G$, we obtain that $E=\leb(G-\bar G)$ and by \eqref{form2}, \eqref{Gbar}, and
\eqref{form2bis} we conclude that
\begin{equation}\label{FEL}
E = \leb\Big(x_2 Ae_2+ x_3 Ae_3\left|\, \partial_2\beta \, \right| \partial_3\beta\Big).
\end{equation}
Equation \eqref{divfree} and Lemma~\ref{lm:EulerQ} guarantee that $\beta(x_1,\cdot)$
is a solution to the problem \eqref{minpb} defining $Q_1(A(x_1))$, for a.e.\ $x_1\in(0,L)$.

As for the first moments, let $\varphi\in C^\infty([0,L];\R^3)$ be such that $\varphi(0)=0$.
Using $\psi(x):=x_2\varphi(x_1)$ as test functions in (\ref{EL2}), we obtain
$$
\int_\Omega ( x_2R^{(h)}E^{(h)}e_1{\,\cdot\,}\varphi'+\tfrac1h R^{(h)}E^{(h)}e_2{\,\cdot\,}\varphi -hx_2
g{\,\cdot\,}\varphi)\, dx =0.
$$
Integrating first with respect to $x_2,x_3$ and using \eqref{section}, the equation reduces to
\begin{equation}\label{el5bis}
\int_0^L (R^{(h)} \tilde E^{(h)}e_1{\,\cdot\,}\varphi' 
+\tfrac1h R^{(h)} \bar E^{(h)}e_2{\,\cdot\,}\varphi)\, dx_1 =0.
\end{equation}
In particular, if we choose $\varphi$ of the form $\varphi=\phi\,R^{(h)} e_1$
with $\phi\in C^\infty([0,L])$ and $\phi(0)=0$, we obtain
\begin{equation}\label{EL5}
\int_0^L ( \phi'\,\tilde E^{(h)}_{11} + \phi\, \tilde E^{(h)}e_1{\,\cdot\,} A^{(h)}e_1
+\phi\, \tfrac1h  \bar E^{(h)}_{12})\, dx_1 =0.
\end{equation}

{}From the estimate (\ref{cons3}) and the identity (\ref{mom0}) it follows that 
the term $\frac1h \bar E^{(h)}_{12}$ is bounded in $L^1(0,L)$. 
Since $A^{(h)}$ and $\tilde E^{(h)}$ are bounded in $L^2((0,L);\mthree)$, 
the product $\tilde E^{(h)}e_1{\,\cdot\,} A^{(h)}e_1$ is also bounded in $L^1(0,L)$.
Therefore, equation (\ref{EL5}) implies that
\begin{equation}\label{EL6}
\|\partial_1\tilde E^{(h)}_{11}\|_{L^1}\le C, \quad \tilde E^{(h)}_{11}(L)=0,
\end{equation}
hence the sequence $\tilde E^{(h)}_{11}$ is strongly compact in $L^p(0,L)$ for every $p<\infty$.

Analogously, one can show that 
\begin{equation}\label{el5ter}
\int_0^L ( R^{(h)} \hat E^{(h)}e_1{\,\cdot\,}\varphi' 
+\tfrac1h R^{(h)} \bar E^{(h)}e_3{\,\cdot\,}\varphi )\, dx_1 =0
\end{equation}
for every $\varphi\in C^\infty([0,L];\R^3)$ such that $\varphi(0)=0$.
Choosing the test function $\varphi$ of the form $\varphi=\phi\,R^{(h)} e_1$
with $\phi\in C^\infty([0,L])$ and $\phi(0)=0$,  one obtains
\begin{equation}\label{EL6.00}
\int_0^L ( \phi'\hat E^{(h)}_{11} - \phi\, \hat E^{(h)}e_1{\,\cdot\,} A^{(h)}e_1
+\phi\, \tfrac1h  \bar E^{(h)}_{13})\, dx_1 =0
\end{equation}
for every $\phi\in C^\infty([0,L])$ with $\phi(0)=0$.
{}From this equation one can deduce as before that
\begin{equation}\label{EL6.1}
\|\partial_1\hat E^{(h)}_{11}\|_{L^1}\le C, \quad \hat E^{(h)}_{11}(L)=0,
\end{equation}
hence the sequence $\hat E^{(h)}_{11}$ is strongly compact in $L^p(0,L)$ for every $p<\infty$.

Finally, let us consider $\phi\, R^{(h)}e_2$ and $\phi\, R^{(h)}e_3$ 
as test functions in \eqref{el5ter} and \eqref{el5bis}, respectively, 
with $\phi\in C^\infty([0,L])$ with $\phi(0)=0$.
Taking the difference of the two equations we obtain
\begin{equation}\label{el5mid}
\begin{array}{c}
\displaystyle
\int_0^L \phi'\,(\hat E^{(h)}_{21}- \tilde E^{(h)}_{31})\, dx_1
- \int_0^L \phi\,(\hat E^{(h)}e_1{\,\cdot\,}A^{(h)}e_2 - \tilde E^{(h)}e_1{\,\cdot\,}A^{(h)}e_3)\, dx_1
\smallskip
\\
\displaystyle
+ \int_0^L \phi\,\tfrac1h (\bar E^{(h)}_{23}-  \bar E^{(h)}_{32})\, dx_1 =0.
\end{array}
\end{equation}
As $A^{(h)}$ and $E^{(h)}$ are bounded in $L^2(\Omega;\mthree)$, the term
$(A^{(h)}\hat E^{(h)})_{21} - (A^{(h)}\tilde E^{(h)})_{31} $ is bounded in $L^1(0,L)$.
The difference $\tfrac1h (\bar E^{(h)}_{23}-  \bar E^{(h)}_{32})$ is also bounded in $L^1(0,L)$
by \eqref{cons3}. Therefore, we deduce from equation \eqref{el5mid} that 
\begin{equation}\label{EL6.2}
\|\partial_1(\hat E^{(h)}_{21}- \tilde E^{(h)}_{31})\|_{L^1}\le C, \quad 
\hat E^{(h)}_{21}(L)- \tilde E^{(h)}_{31}(L)=0,
\end{equation}
hence the sequence $\hat E^{(h)}_{21}- \tilde E^{(h)}_{31}$
is strongly compact in $L^p(0,L)$ for every $p<\infty$.
\bigskip

\noindent
{\em Step 6. Convergence of the energy by the div-curl lemma.}
\smallskip

\noindent
The strong compactness of the sequences $(\tilde E^{(h)}_{11})$, $(\hat E^{(h)}_{11})$, and
$(\hat E^{(h)}_{21}-\hat E^{(h)}_{31})$ allows us to pass to the limit in the energy integral
$$
\frac{1}{h^2}\int_\Omega DW(Id+hG^{(h)}){\,:\,}hG^{(h)}\, dx = 
\int_\Omega E^{(h)}{\,:\,}G^{(h)}\, dx.
$$
This can be done by exploiting the div-curl structure of the product $E^{(h)}{\,:\,}G^{(h)}$; indeed, 
the Euler-Lagrange equation (\ref{EL2}) asserts that
the scaled divergence of $R^{(h)}E^{(h)}$ is infinitesimal in $L^2(\Omega;\R^3)$ as $h\to 0$,
while the decomposition (\ref{form1}) guarantees that the matrix
$R^{(h)}G^{(h)}$ has basically the structure of a scaled gradient.

Let us fix $\varphi\in C^\infty(0,L)$ with $\varphi(0)=0$. Using formula (\ref{form1}) we have
\begin{equation}\label{div1}
\begin{array}{c}
\displaystyle
\int_\Omega \varphi E^{(h)}{\,:\,}G^{(h)}\, dx = 
\int_\Omega \varphi R^{(h)}E^{(h)}{\,:\,}R^{(h)}G^{(h)}\, dx 
\smallskip 
\\
\displaystyle
=  \int_\Omega \varphi R^{(h)}E^{(h)}{\,:\,}\nabla_h z^{(h)}\, dx
+ \int_\Omega \varphi E^{(h)}e_1{\,\cdot\,}(x_2A^{(h)}e_2+x_3A^{(h)}e_3) \, dx. 
\end{array}
\end{equation}
Concerning the first term on the right-hand side, the Euler-Lagrange equation (\ref{EL2}) yields
$$
\int_\Omega \varphi R^{(h)}E^{(h)}{\,:\,}\nabla_h z^{(h)}\, dx =
h\int_\Omega \varphi g{\,\cdot\,} z^{(h)}\, dx 
- \int_\Omega \varphi'R^{(h)}E^{(h)}e_1{\,\cdot\,} z^{(h)}\, dx.
$$
Since $z^{(h)}\to z$ strongly in $L^2(\Omega; \R^3)$ and $R^{(h)}E^{(h)}\wto RE$
weakly in $L^2(\Omega;\mthree)$, we can pass to the limit
in the formula above and we get
$$
\lim_{h\to 0}\int_\Omega \varphi R^{(h)}E^{(h)}{\,:\,}\nabla_h z^{(h)}\, dx = 
- \int_\Omega \varphi'RE e_1{\,\cdot\,} z\, dx.
$$
Taking into account the fact that $z$ is independent of $x_2,x_3$ and using
the identity (\ref{mom1}), we have
$$
\int_\Omega \varphi'RE e_1{\,\cdot\,} z\, dx =
\int_0^L \varphi'R\bar E e_1{\,\cdot\,} z\, dx_1= 0, 
$$
hence
\begin{equation} \label{div2}
\lim_{h\to 0}\int_\Omega \varphi R^{(h)}E^{(h)}{\,:\,}\nabla_h z^{(h)}\, dx = 0.
\end{equation}

As for the last term in (\ref{div1}), integrating first with respect to the cross-section variables
we have
$$
\begin{array}{c}
\displaystyle
\int_\Omega \varphi E^{(h)}e_1{\,\cdot\,}(x_2A^{(h)}e_2+x_3A^{(h)}e_3) \, dx 
\smallskip
\\
\displaystyle
= \int_0^L \varphi ( \tilde E^{(h)}_{11}A^{(h)}_{12} +
\hat E^{(h)}_{11}A^{(h)}_{13} ) \, dx_1 +
 \int_0^L \varphi (\hat E^{(h)}_{21} -\tilde E^{(h)}_{31})A^{(h)}_{23}\, dx_1.
\end{array}
$$
As $\tilde E^{(h)}_{11}$, $\hat E^{(h)}_{11}$, and 
$\hat E^{(h)}_{21}-\tilde E^{(h)}_{31}$ are strongly compact in $L^2(0,L)$ by Step~5,
we can pass to the limit and we obtain
\begin{equation}\label{div3}
\begin{array}{c}
\displaystyle
\lim_{h\to 0}\int_\Omega \varphi E^{(h)}e_1{\,\cdot\,}(x_2A^{(h)}e_2+x_3A^{(h)}e_3) \, dx
\smallskip
\\
\displaystyle 
= \int_0^L \varphi ( \tilde E_{11}A_{12} + \hat E_{11}A_{13} )\, dx_1 +
\int_0^L \varphi(\hat E_{21} -\tilde E_{31})A_{23}\, dx_1
\smallskip
\\
\displaystyle 
=\int_\Omega \varphi Ee_1{\,\cdot\,}(x_2Ae_2+x_3Ae_3) \, dx.
\end{array}
\end{equation}
Now from the first equality in (\ref{form2}) it follows that
\begin{equation}\label{div4.1}
\int_\Omega \varphi Ee_1{\,\cdot\,}(x_2Ae_2+x_3Ae_3) \, dx
= \int_\Omega \varphi Ee_1{\,\cdot\,}Ge_1\, dx
- \int_\Omega \varphi Ee_1{\,\cdot\,}R^Tz'\, dx.
\end{equation}
Since $R^Tz'$ does not depend on $x_2,x_3$, identity (\ref{mom1}) implies
$$
\int_\Omega \varphi Ee_1{\,\cdot\,}  R^Tz'\, dx =
 \int_0^L \varphi \bar Ee_1{\,\cdot\,}  R^Tz'\, dx_1 =0.
$$
Thus, equality \eqref{div4.1} reduces to
\begin{equation}\label{div4}
\int_\Omega \varphi Ee_1{\,\cdot\,}(x_2Ae_2+x_3Ae_3) \, dx
= \int_\Omega \varphi Ee_1{\,\cdot\,}Ge_1\, dx.
\end{equation}
Combining together (\ref{div1})--(\ref{div3}) and (\ref{div4}), 
we conclude that 
\begin{equation}\label{en-conv0}
\lim_{h\to 0}\int_\Omega \varphi E^{(h)}{\,:\,}G^{(h)}\, dx
= \int_\Omega \varphi Ee_1 {\,\cdot\,}Ge_1 \, dx
\end{equation}
By \eqref{divfree} the matrix $(Ee_2\,|\,Ee_3)$ is divergence free in $S$ with zero normal component
on $\partial S$ for a.e.\ $x_1\in(0,L)$, while $(Ge_2\,|\,Ge_3)$ is a gradient
by \eqref{form2bis}. As the test function $\varphi$ depends only on the variable $x_1$,
the divergence theorem yields
$$
\int_\Omega \varphi ( Ee_2{\,\cdot\,}Ge_2 + Ee_3{\,\cdot\,}Ge_3 )\, dx=0,
$$
hence
\begin{equation}\label{div45}
\int_\Omega \varphi E{\,:\,}G\, dx =\int_\Omega \varphi Ee_1{\,\cdot\,}Ge_1 \, dx.
\end{equation}
By (\ref{en-conv0}) and (\ref{div45}) 
we finally obtain the convergence of the energies
\begin{equation}\label{en-conv}
\lim_{h\to 0}\int_\Omega \varphi E^{(h)}{\,:\,}G^{(h)}\, dx
= \int_\Omega \varphi E {\,:\,}G \, dx
\end{equation}
for every $\varphi\in C^\infty_0(0,L)$.
\bigskip

\noindent
{\em Step 7. Definition of the truncated deformations.}
\smallskip

\noindent
In order to pass to the limit in the Euler-Lagrange equations \eqref{EL5}, \eqref{EL6.00}, 
and \eqref{el5mid}, a strong $L^2$-compactness for the sequence $(E^{(h)})$ is required.
If $hG^{(h)}$ converges to $0$ uniformly, then by Taylor expansion one can replace
$E^{(h)}$ by $\leb G^{(h)}$ in \eqref{en-conv}. Using the fact that $E=\leb G$ and $\leb$ is
positive definite on symmetric matrices, one can conclude strong convergence for $\sym G^{(h)}$
and hence of $E^{(h)}$, outside a neighbourhood of $x_1=0$ (see Step~7 of the proof
of Theorem~1.1 in \cite{MMS06}). 

To avoid the extra assumption $h\|G^{(h)}\|_\infty\to 0$, we introduce an auxiliary sequence
of truncated deformations $u^{(h)}$, whose corresponding scaled strains $H^{(h)}$
satisfy $h\|H^{(h)}\|_\infty\to 0$ (see \eqref{tildeG}).
The main point will be then to show strong convergence of $\sym H^{(h)}$ (see Step~8).
This will imply, as before, strong convergence of the corresponding truncated stress $F^{(h)}$
(outside a neighbourhood of $x_1=0$). To pass to the limit in the Euler-Lagrange equations
and conclude the proof, we will then need to estimate the remainder term $E^{(h)}-F^{(h)}$.
This will be done in Step~9, using again the div-curl lemma and exploiting our careful choice of the truncations.

To carry out this plan, we consider the functions $z^{(h)}$ defined in (\ref{def:z^h}) and their rescalings
$\check z^{(h)}(x):=z^{(h)}(x_1,\frac{x_2}{h},\frac{x_3}{h})$.
Applying Lemma~\ref{truncate} to $\check z^{(h)}$ with $a=h^{-5/8}$ and $b=h^{-7/8}$
and undoing the rescaling, we construct a new sequence of functions
$w^{(h)}:\Omega\to\R^3$ with the following properties:
\begin{equation}\label{tilde1}
\|\nabla_h w^{(h)}\|_{L^\infty}\le \lambda_h, 
\end{equation}
\begin{eqnarray}
\lambda_h^2\leb^3(N_h) & \leq &
\frac{C}{\ln(1/h)}\int_\Omega |\nabla_h z^{(h)}|^2\, dx 
\nonumber
\\
& \leq &
\frac{C}{\ln(1/h)}
\int_{\Omega\vphantom{_O}} (|G^{(h)}|^2+|A^{(h)}|^2 )\, dx, \label{tilde2}
\\
\|\nabla_h z^{(h)} - \nabla_h w^{(h)} \|^2_{L^2} & \leq & 
\frac{C}{\ln(1/h)}\int_\Omega |\nabla_h z^{(h)}|^2\, dx,
\label{tilde3}
\end{eqnarray}
where $\lambda_h\in[h^{-5/8},h^{-7/8}]$ and 
$N_h :=\{x\in \Omega: \ z^{(h)}(x)\neq w^{(h)}(x)\}$.
In particular we have
\begin{equation}\label{tilde4}
h^{1/2}\lambda_h\to \infty, \quad h\lambda_h\to 0, \quad \hbox{and} \quad
\lambda_h^2  \leb^3(N_h) \to 0.
\end{equation}

We can introduce now the sequence of approximated deformations $u^{(h)}:\Omega\to\R^3$, which
are associated with the auxiliary functions $w^{(h)}$:
$$
u^{(h)}:=hw^{(h)}+\int_0^{x_1}R^{(h)}(s)e_1\, ds+hx_2R^{(h)}e_2+hx_3R^{(h)}e_3.
$$
Let $H^{(h)}:\Omega\to\mthree$ be the corresponding approximated strains defined by
the relation
$$
\nabla_hu^{(h)} = R^{(h)} (Id + h  H^{(h)}),
$$
and let $F^{(h)}:\Omega\to\mthree$ be the corresponding stresses defined as
\begin{equation}\label{defEtilde}
F^{(h)}:= \frac1h DW(Id + h  H^{(h)}).
\end{equation}
Using the definition of $u^{(h)}$ it is easy to see that
\begin{equation}\label{defGtilde}
H^{(h)} = (R^{(h)})^T\nabla_hw^{(h)}+x_2 A^{(h)} e_2\otimes e_1
+ x_3 A^{(h)} e_3\otimes e_1.
\end{equation}
It follows from \eqref{tilde3} that $\nabla_h w^{(h)}$ and $\nabla_h z^{(h)}$
have the same weak limit and hence by \eqref{nablaz2} 
\begin{equation}\label{Hweak}
H^{(h)}\wto G \quad \hbox{weakly in } L^2(\Omega;\mthree).
\end{equation}
\medskip

\noindent
{\em Step 8. $L^\infty$-convergence of $hH^{(h)}$ and strong convergence of $\sym H^{(h)}$ and $F^{(h)}$.}
\smallskip

\noindent
We recall the estimate
\begin{equation}\label{stima}
\sup|f-\bar f|^2\le 2\| f\|_{L^2} \|f'\|_{L^2}, \quad \text{with } \bar f:=\frac1L \int_0^L f\, dx.
\end{equation}
As $(R^{(h)})'$ and $h(R^{(h)})''$ are bounded in $L^2(0,L)$ by (\ref{rig2}), we deduce that
$|(R^{(h)})' |\le Ch^{-1/2}$, and therefore $|A^{(h)}|\le Ch^{-1/2}$. This inequality and (\ref{tilde1}) imply that
\begin{equation}\label{tildeG}
h|H^{(h)}|\le Ch\lambda_h +Ch^{1/2}\to 0.
\end{equation}
By Taylor expansion of $DW$ around the identity matrix we have
\begin{equation}\label{Eexp0}
F^{(h)}= \frac1h DW(Id+hH^{(h)})= \leb H^{(h)} +\frac1h\eta(hH^{(h)}),
\end{equation}
where $|\eta(A)|/|A|\to 0$, as $|A|\to 0$. For every $t>0$ let us define
$$
\omega(t):=\sup\Big\{\frac{|\eta(A)|}{|A|}: \ |A|\le t \Big\};
$$
then, it is easy to see that $\omega(t)\to 0$, as $t\to 0^+$. 
The expansion \eqref{Eexp0} and the definition of $\omega$ yield
$$
|\leb H^{(h)}{\,:\,}H^{(h)} - F^{(h)} {\,:\,} H^{(h)}|
\leq \omega(h\|H^{(h)}\|_{L^\infty})|H^{(h)}|^2.
$$
Together with \eqref{tildeG}, we obtain for every $\varphi\in C^\infty([0,L])$ 
\begin{equation}\label{LH}
\int_\Omega \varphi \leb H^{(h)}{\,:\,}H^{(h)}\, dx
- \int_\Omega \varphi  F^{(h)} {\,:\,} H^{(h)}\, dx
\to 0.
\end{equation}

We now claim that
\begin{equation}\label{FhHh}
\int_\Omega \varphi  F^{(h)} {\,:\,} H^{(h)}\, dx - \int_\Omega \varphi  E^{(h)} {\,:\,} G^{(h)}\, dx
\to 0.
\end{equation}
Combining together the convergence of energy \eqref{en-conv}, the weak convergence \eqref{Hweak},
and \eqref{LH}, this would imply that
\begin{equation}\label{GhHh}
\lim_{h\to 0} \int_\Omega \varphi \leb(H^{(h)}-G){\,:\,}(H^{(h)}-G)\, dx=0
\end{equation}
for every $\varphi\in C^\infty([0,L])$ with $\varphi(0)=0$.
{}From the assumptions on $W$ we infer that there exists a constant $C>0$ such that
$$
\leb A{\,:\,}A\geq C|\sym A|^2
$$
for every $A\in\mthree$. This inequality, together with \eqref{GhHh}, implies
\begin{equation}\label{conv-tildeG}
\sym (H^{(h)}-G) \to 0 \qquad \hbox{strongly in } L^2((a,L){\times}S;\mthree)
\end{equation}
for every $a>0$. Using again the Taylor expansion \eqref{Eexp0}, we easily deduce that
\begin{equation}\label{convH}
F^{(h)} \to E \qquad \hbox{strongly in } L^2((a,L){\times}S;\mthree).
\end{equation}

In order to prove \eqref{FhHh} we write the difference as
\begin{equation}\label{diff}
\int_\Omega \varphi E^{(h)} {\,:\,} (H^{(h)}-G^{(h)})\, dx
+ \int_{\Omega} \varphi  ( F^{(h)} -E^{(h)} ){\,:\,}H^{(h)}\, dx.
\end{equation}
The first term can be controlled by the div-curl lemma; indeed, equalities (\ref{defGtilde}) and (\ref{form1}) yield
$$
R^{(h)}(H^{(h)}- G^{(h)})= \nabla_h(w^{(h)} -z^{(h)}),
$$
so that, by the Euler-Lagrange equation (\ref{EL2}), we have
$$
\begin{array}{c}
\displaystyle
\int_{\Omega} \varphi E^{(h)}
{\,:\,}(H^{(h)}- G^{(h)})\, dx
=  \int_{\Omega} \varphi R^{(h)}E^{(h)}
{\,:\,}\nabla_h(w^{(h)} - z^{(h)})\, dx
\smallskip
\\
\displaystyle
= h\int_{\Omega}\varphi g{\,\cdot\,} (w^{(h)} - z^{(h)})\, dx
-\int_{\Omega} \varphi'R^{(h)}E^{(h)}e_1{\,\cdot\,} (w^{(h)} - z^{(h)})\, dx.
\end{array}
$$
Since the sequence $w^{(h)} - z^{(h)}$ converges to $0$ strongly in $L^2(\Omega;\R^3)$ and
$R^{(h)}E^{(h)}$ is bounded in $L^2(\Omega;\mthree)$, we conclude that
$$
\lim_{h\to 0} 
\int_{\Omega} \varphi E^{(h)}{\,:\,}(H^{(h)}- G^{(h)})\, dx
= 0.
$$
To estimate the second integral in \eqref{diff} we recall that $F^{(h)}$ and $E^{(h)}$
are bounded in $L^2(\Omega;\mthree)$. Therefore, by H\"older inequality and by \eqref{tildeG}
we have 
$$
\int_{\Omega} |\varphi (F^{(h)}- E^{(h)})
{\,:\,}H^{(h)}| \, dx \le C\Big( \int_{N_h} |H^{(h)}|^2\, dx \Big)^{1/2}
\leq C[(\lambda_h^2+h^{-1})\leb^3(N_h)]^{1/2}.
$$
As the right-hand side converges to zero by \eqref{tilde4},
this concludes the proof of the claim \eqref{FhHh}.
\bigskip

\noindent
{\em Step 9. Passage to the limit in the Euler-Lagrange equations.}
\smallskip

\noindent
Let us fix $\phi\in C^\infty([0,L])$ vanishing on an interval $(0,a)$. In order to pass to the
limit in the Euler-Lagrange equations, we need to prove some preliminary convergence results.
First of all we claim that
\begin{equation}\label{a1}
\lim_{h\to 0}\int_\Omega \phi\, x_k E^{(h)}e_1{\,\cdot\,}A^{(h)}e_j \, dx
= \int_\Omega \phi\, x_k Ee_1{\,\cdot\,}Ae_j \, dx
\end{equation}
for every $k=2,3$ and every $j=1,2,3$.
Indeed, 
\begin{eqnarray}
\int_\Omega \phi \, x_k E^{(h)}e_1{\,\cdot\,}A^{(h)}e_j \, dx
& = & \int_\Omega \phi \, x_k F^{(h)}e_1{\,\cdot\,}A^{(h)}e_j \, dx 
\nonumber \\
& &
+ \int_{N_h} \phi\, x_k (E^{(h)}-F^{(h)})e_1{\,\cdot\,}A^{(h)}e_j \, dx.
\label{a02}
\end{eqnarray}
By the strong convergence (\ref{convH}) we have
$$
\int_\Omega \phi\, x_k F^{(h)}e_1{\,\cdot\,}A^{(h)}e_j \, dx \to  \int_\Omega 
\phi\, x_k Ee_1{\,\cdot\,}Ae_j \, dx.
$$
As for the last term in \eqref{a02}, using H\"older's inequality we obtain
$$
\int_{N_h} |\phi\, x_k (E^{(h)}-F^{(h)})e_1{\,\cdot\,}A^{(h)}e_j| \, dx \le
C\Big( \int_{N_h} |A^{(h)} |^2\, dx  \Big)^{1/2}.
$$
Since $|A^{(h)}|\le Ch^{-1/2}$ and $h^{-1}\leb^3(N_h)\to 0$ by (\ref{tilde4}),
the previous estimate implies that the second integral on the right-hand side of (\ref{a02})
converges to $0$. This concludes the proof of the claim \eqref{a1}.

Integrating first with respect to the variables of the cross-section in
\eqref{a1}, we obtain for $k=2,3$ and $j=1$
\begin{eqnarray}
& \displaystyle\label{a1tilde}
\lim_{h\to 0}\int_0^L \phi\, \tilde E^{(h)}e_1{\,\cdot\,}A^{(h)}e_1 \, dx_1 = 
\int_0^L \phi\, \tilde Ee_1{\,\cdot\,}Ae_1 \, dx_1,
\\
& \displaystyle \label{a1hat}
\lim_{h\to 0}\int_0^L  \phi\, \hat E^{(h)}e_1{\,\cdot\,}A^{(h)}e_1 \, dx_1 = 
\int_0^L  \phi\, \hat Ee_1{\,\cdot\,}Ae_1 \, dx_1.
\end{eqnarray}

Arguing as in the proof of \eqref{a1}, it is easy to show that 
\begin{equation}\label{otim}
\lim_{h\to 0}\int_0^L \phi\,x_k\, \skw ( E^{(h)}e_1\otimes A^{(h)}e_k)\, dx =
\int_\Omega \phi\, x_k\, \skw ( Ee_1\otimes Ae_k)\, dx
\end{equation}
for every $k=2,3$ and every $\phi\in C^\infty([0,L])$ vanishing on $(0,a)$.

In order to pass to the limit in the Euler-Lagrange equations (\ref{EL5})
and \eqref{EL6.00}, it remains to study the convergence of the terms
$$
\int_0^L \phi\,\tfrac1h \bar E^{(h)}_{1k}\, dx_1
$$
for $k=2,3$. We first decompose the integral as 
\begin{equation}\label{pass}
\int_0^L \phi\,\tfrac1h \bar E^{(h)}_{1k}\, dx_1=
\int_0^L \phi\,\tfrac1h \bar E^{(h)}_{k1}\, dx_1 +
2\int_\Omega \phi\, \tfrac1h\skw(E^{(h)})_{1k}\, dx.
\end{equation}
By \eqref{mom0} we immediately deduce that
\begin{equation}\label{noto}
\lim_{h\to 0}\int_0^L \phi\,\tfrac1h \bar E^{(h)}_{k1}\, dx_1 =
-\int_0^L \phi\, R^T\tilde g{\,\cdot\,}e_k\, dx_1.
\end{equation}
As for the second integral on the right-hand side of \eqref{pass},
it follows from (\ref{symm}) that $\tfrac1h \skw E^{(h)} = - \skw (E^{(h)}(G^{(h)})^T)$, 
so that equality (\ref{form1}) yields
\begin{eqnarray*}
\tfrac1h \skw(E^{(h)}) & = & - \skw ( E^{(h)} (\nabla_h z^{(h)})^T R^{(h)} )
 - x_2\, \skw ( E^{(h)} e_1\otimes A^{(h)}e_2) 
\\
& &
- x_3\, \skw ( E^{(h)} e_1\otimes A^{(h)}e_3).
\end{eqnarray*}
Since $\skw A=\skw(RAR^T)$ for every $A\in\mthree$ and every $R\in SO(3)$, 
we have that
$$
\skw ( E^{(h)} (\nabla_h z^{(h)})^T R^{(h)} )=\skw ( R^{(h)}E^{(h)} (\nabla_h z^{(h)})^T).
$$
This identity, together with the Euler-Lagrange equation \eqref{EL2} and the strong convergence of
$z^{(h)}$, implies that
\begin{eqnarray*}
\lefteqn{\lim_{h\to 0}\int_\Omega \skw ( E^{(h)} (\nabla_h z^{(h)})^T R^{(h)} ) \phi\, dx}\\
& = & \int_\Omega \skw ( REe_1\otimes z ) \phi' \, dx 
= \int_0^L \skw ( R\bar E e_1\otimes z ) \phi' \, dx =0,
\end{eqnarray*}
where we have used the fact that $z$ and $R$ are independent of $x_2,x_3$ and that $\bar Ee_1=0$
by \eqref{mom1}.
Combining this equality with \eqref{otim}, we conclude that
\begin{equation}\label{b1}
\lim_{h\to 0}\int_0^L \phi\, \tfrac1h \skw \bar E^{(h)} \, dx_1= 
-\int_0^L \phi\, \skw(\tilde Ee_1\otimes Ae_2+ \hat Ee_1\otimes Ae_3) \, dx_1 
\end{equation}
for every $\phi\in C^\infty([0,L])$ vanishing on $(0,a)$.

By \eqref{noto} and \eqref{b1} we finally obtain that
\begin{eqnarray*}
\lefteqn{\lim_{h\to 0}\int_0^L \phi\, \tfrac1h \bar E^{(h)}_{1k} \, dx_1}
\\
& = & \hspace{-2mm}
-2\int_0^L \phi\, \skw(\tilde Ee_1\otimes Ae_2+ \hat Ee_1\otimes Ae_3)_{1k} \, dx_1 
 -\int_0^L \phi\, R^T\tilde g{\,\cdot\,} e_k \, dx_1.
\end{eqnarray*}
Together with \eqref{a1tilde} and \eqref{a1hat}, this shows that we can pass to the limit
in (\ref{EL5}) and \eqref{EL6.00}. Thus, we obtain the equations 
\begin{equation}\label{ELfin1}
\int_0^L ( \phi'\, \tilde E_{11} + \phi\, A_{13}(\hat E_{21} -\tilde E_{31}) - \phi\, A_{23}\hat E_{11}
- \phi\, R^T\tilde g{\,\cdot\,} e_2)\, dx_1 =0
\end{equation}
and
\begin{equation}\label{ELfin2}
\int_0^L ( \phi'\, \hat E_{11} - \phi\, A_{12}(\hat E_{21} -\tilde E_{31}) + \phi\, A_{23}\tilde E_{11}
- \phi\, R^T\tilde g{\,\cdot\,} e_3)\, dx_1 =0
\end{equation}
for every $\phi\in C^\infty([0,L])$ vanishing on $(0,a)$.

Analogously, by \eqref{a1} we deduce
\begin{eqnarray*}
\lefteqn{\lim_{h\to 0}\int_0^L \phi\, (\hat E^{(h)}e_1{\,\cdot\,} A^{(h)}e_2 
- \tilde E^{(h)}e_1{\,\cdot\,} A^{(h)}e_3)\, dx_1}
\\
& = & \int_0^L \phi\, (\hat Ee_1{\,\cdot\,} Ae_2 - \tilde Ee_1{\,\cdot\,} Ae_3) \, dx_1,
\end{eqnarray*}
while by \eqref{b1} we have
$$
\lim_{h\to 0}\int_0^L \phi\, \tfrac1h (\bar E^{(h)}_{23}-  \bar E^{(h)}_{32}) \, dx_1
= - \int_0^L \phi\, (A_{32} \tilde E_{21} - A_{23} \hat E_{31}) \, dx_1.
$$
Combining these two properties,
we can pass to the limit also in the equation \eqref{el5mid} and we obtain
\begin{equation}\label{ELfin3}
\int_0^L (\phi'\, (\hat E_{21}- \tilde E_{31})+ \phi\, A_{12} \hat E_{11} - \phi\, A_{13}\tilde E_{11}) 
\, dx_1=0
\end{equation}
for every $\phi\in C^\infty([0,L])$ vanishing on $(0,a)$.

By approximation it is easy to see that the limiting equations \eqref{ELfin1}, \eqref{ELfin2},
and \eqref{ELfin3} hold for every $\phi\in C^\infty([0,L])$ with $\phi(0)=0$.

Finally, taking into account \eqref{FEL} and integrating by parts, 
one can check that conditions \eqref{ELfin1}--\eqref{ELfin3}
coincide with the Euler-Lagrange equations \eqref{eq:EL2} for $J_2$.
\qed

\section{Truncation and compactness}
\label{sec:trcp}

In this section we collect some auxiliary results which were used in the proof of Theorem~\ref{thm:1.1}.

The first proposition contains an approximation result  by means of smooth rotations
for sequences of deformations with elastic energy of order $h^2$. 
This is the point where the rigidity lemma by Friesecke, James, and 
M\"uller (see \cite[Theorem~3.1]{FJM02}) is used in a crucial way.

\begin{proposition}\label{rotation}
Let $(u^{(h)})\subset W^{1,2}(\Omega;\R^3)$ be a sequence such that
$$
F^{(h)}(u^{(h)}) := \int_\Omega \dist^2(\nabla_h u^{(h)}, SO(3))\, dx\le Ch^2
$$
for every $h>0$. 
Then there exists an associated sequence $(R^{(h)})\subset C^\infty((0,L);\mthree)$ such that 
\begin{eqnarray}
& \displaystyle 
R^{(h)} (x_1)\in SO(3) \quad \hbox{for every }x_1\in(0,L), \label{rot1}
\\
& \displaystyle 
\|\nabla_h u^{(h)} - R^{(h)}\|_{L^2} \le Ch, \label{rot2}
\\
& \displaystyle 
\|(R^{(h)})'\|_{L^2}+
h\, \| (R^{(h)})''\|_{L^2}\le  C  \label{rot3} 
\end{eqnarray}
for every $h>0$. If, in addition, 
$u^{(h)}(0,x_2,x_3) = (0, hx_2,hx_3)$, then
\begin{equation} \label{bcrigid}
| R^{(h)}(0) - Id | \leq C \sqrt{h}.
\end{equation}
\end{proposition}

\begin{proof}
The argument follows closely the proof of \cite[Proposition~4.1]{MMS06}.
For every $h>0$ the set $\Omega_h$ can be partitioned in cylinders of the form $I_h{\times}hS$,
where $I_h$ is an interval of length comparable to $h$.
Applying the rigidity estimate \cite[Theorem~3.1]{FJM02} in each such cylinder,
one first construct a sequence $(Q^{(h)})$ of piecewise constant rotations satisfying \eqref{rot2}
and a difference quotient variant of \eqref{rot3}.
As the mollifications $\tilde Q^{(h)}$ of $Q^{(h)}$ at scale $h$ are uniformly close to $Q^{(h)}$,
it is possible to project $\tilde Q^{(h)}$ back on $SO(3)$; this provides the sequence $(R^{(h)})$.
For the details we refer to~\cite{MMS06}.
\end{proof}

The next proposition allows to identify the weak limit of the sequence of stresses $(E^{(h)})$,
once the weak limit of the strains $(G^{(h)})$ is known. For the proof, which is based on Taylor
expansion, we refer to \cite[Proposition~4.2]{MMS06}.

\begin{proposition} \label{Ehweak}
Assume that the energy density $W$ is differentiable and 
its derivative $DW$ is Lipschitz continuous. Assume moreover
that $DW$ is differentiable at the identity.
Suppose that 
$$
G^{(h)} \wto G \quad \mbox{weakly in } L^2(\Omega;\mthree)
$$
and define the rescaled stresses as in  (\ref{Ehdef}) by
$$
E^{(h)} := \frac1h DW(Id + h G^{(h)}). 
$$
Then
\begin{equation}
E^{(h)}\wto E:=\leb\, G \quad \mbox{weakly in }L^2(\Omega;\mthree),
\end{equation}
where $\leb:=D^2W(Id)$.
\end{proposition}

We conclude this section with the truncation lemma used in the proof of Theorem~\ref{thm:1.1}.
This a variant for thin domains of the standard results on the truncations of gradients
(see, e.g., \cite{EG}).

\begin{lemma}\label{truncate}
There exists a constant $C>0$ with the following property: 
for every $h>0$, every $b>a>0$ and every $u\in W^{1,2}(\Omega_h;\R^3)$
there exist $\lambda\in [a,b]$ and a function $v\in W^{1,\infty}(\Omega_h;\R^3)$ such that
\begin{eqnarray}
& \displaystyle
\label{trunc1}
\|\nabla v\|_{L^\infty}\le \lambda, \vphantom{\int}
\\
& \displaystyle
\label{trunc2}
\lambda^2\leb^3(\{x\in \Omega_h: \ u(x)\neq v(x)\})\le
\frac{C}{\ln(b/a)}\int_{\{x\in \Omega_h:\ |\nabla u(x)|>\lambda\}} |\nabla u|^2\, dx,
\\
& \displaystyle
\label{trunc3}
\|\nabla u-\nabla v\|^2_{L^2}\le 
\frac{C}{\ln(b/a)}\int_{\{x\in \Omega_h:\ |\nabla u(x)|>\lambda\}} |\nabla u|^2\, dx.
\end{eqnarray}
\end{lemma}

\begin{proof}
Let $Q$ be a square containing $S$. Without loss of generality we can assume that $Q=(0,M)^2$.
Let 
$$
V:=\Big\{ v\in L^2(S;\R^3): \ \bar v:=\int_S v\,dx_2dx_3=0\Big\}.
$$
Then there exists a linear extension operator
$\tilde\E: V\to \{ v\in L^2(Q;\R^3): {\rm supp}\, v\subset\subset Q\}$ such that
$\tilde\E(v)\in W^{1,2}(Q;\R^3)$ for every $v\in V\cap W^{1,2}(S;\R^3)$ and
for some constant $C>0$ there holds
\begin{eqnarray}
& \|\tilde\E(v)\|_{L^2(Q)}\leq C \|v\|_{L^2(S)} \quad \text{for every } v\in V, \label{1ext}
\\
& \|\nabla\!_{x_2,x_3}\tilde\E (v)\|_{L^2(Q)}\leq C \|\nabla\!_{x_2,x_3}v\|_{L^2(S)}
 \quad \text{for every } v\in V\cap W^{1,2}(S;\R^3) \label{2ext}
\end{eqnarray}
(see, e.g., \cite{Stein}). We can extend $\tilde\E$ to the whole space $L^2(S;\R^3)$ by 
considering the operator $\E:L^2(S;\R^3)\to L^2(Q;\R^3)$ defined by
$$
\E(v):=\tilde\E(v-\bar v)+\bar v \quad \text{for every } v\in  L^2(S;\R^3).
$$
It is easy to see that, if $v\in W^{1,2}(S;\R^3)$, then $\E(v)-\bar v\in W^{1,2}_0(Q;\R^3)$.
Moreover, it follows immediately from \eqref{1ext} and \eqref{2ext} that there exists a constant $C$
such that
\begin{eqnarray}
& \|\E(v)\|_{L^2(Q)}\leq C \|v\|_{L^2(S)} \quad \text{for every } v\in L^2(S;\R^3), \label{11ext}
\\
& \|\nabla\!_{x_2,x_3}\E (v)\|_{L^2(Q)}\leq C \|\nabla\!_{x_2,x_3}v\|_{L^2(S)}
 \quad \text{for every } v\in W^{1,2}(S;\R^3) \label{22ext}.
\end{eqnarray}

Let $h>0$ and let $\E_h:L^2(hS;\R^3)\to L^2(hQ;\R^3)$ be the extension operator obtained by scaling $\E$. Then, inequalities \eqref{11ext} and \eqref{22ext} imply that
\begin{eqnarray}
& \|\E_h(v)\|_{L^2(hQ)}\leq C \|v\|_{L^2(hS)} \quad \text{for every } v\in L^2(hS;\R^3), \label{h1ext}
\\
& \|\nabla\!_{x_2,x_3}\E_h (v)\|_{L^2(hQ)}\leq C \|\nabla\!_{x_2,x_3}v\|_{L^2(hS)}
 \quad \text{for every } v\in W^{1,2}(hS;\R^3) \label{h2ext},
\end{eqnarray}
where the constant $C$ is independent of $h$.

Now, let $u\in W^{1,2}(\Omega_h;\R^3)$. First of all we can extend $u$ to the set 
$U_h:=(0,L){\times}hQ$ by defining
$$
\tilde u(x_1,\cdot):=\E_h(u(x_1,\cdot))
$$
for a.e.\ $x_1\in(0,L)$. By \eqref{h1ext} and \eqref{h2ext} we deduce that
\begin{eqnarray}
& \|\tilde u\|_{L^2(U_h)}\leq C \|u\|_{L^2(\Omega_h)}, \label{uh}
\\
& \|\nabla\!_{x_2,x_3}\tilde u \|_{L^2(U_h)}\leq C \|\nabla\!_{x_2,x_3}u\|_{L^2(\Omega_h)}. 
\label{nablau}
\end{eqnarray}
As $\E_h$ is a linear operator, we have that $\partial_1 \tilde u(x_1,\cdot)=\E_h(\partial_1 u(x_1,\cdot))$
for a.e.\ $x_1\in(0,L)$, and thus, by \eqref{h1ext}
\begin{equation}\label{ud1}
\|\partial_1\tilde u \|_{L^2(U_h)}\leq C \|\partial_1 u\|_{L^2(\Omega_h)}.
\end{equation}
As $\tilde u$ is constant on $(0,L){\times}h\partial Q$, we can extend 
$\tilde u$ by successive reflection to the set $U:=(0,L){\times}Q$. By \cite[Lemma~4.3]{MMS06}
there exist $\lambda \in [a,b]$ 
and $w \in W^{1,\infty}(U;\R^3)$ such that
\begin{equation}\label{wOm}
\|\nabla w\|_{L^\infty(U)}\le \lambda
\end{equation}
and 
\begin{equation}\label{wOm2}
\lambda^2 \leb^3(\{x\in U: \ \tilde u(x)\neq w(x)\}) \le 
\frac{C}{\ln(b/a)}
\int_{U} |\nabla \tilde u|^2\, dx.
\end{equation}

Let $N_h$ be the largest integer such that $h(N_h+1)\leq 1$. For $i,j=0,\dots,N_h$ let $Q_{h,ij}$ be the square $(ihM,jhM)+hQ$, let $S_{h,ij}:=(0,L){\times}Q_{h,ij}$, and let 
$$
R_h:= U\setminus  \bigcup_{0\leq i,j\leq N_h} S_{h,ij}.
$$
Since
$$
\sum_{0\leq i,j \leq N_h} \leb^3(\{\tilde u\neq w\}\cap S_{h,ij})\le  \leb^3(\{\tilde u\neq w\}),
$$
there exists some indeces $i_0,j_0$ such that 
\begin{eqnarray}
\lambda^2 \leb^3(\{\tilde u\neq w\}\cap S_{h,i_0j_0})
& \leq & \frac{1}{(N_h+1)^2} \lambda^2 \leb^3(\{\tilde u\neq w\}) \nonumber 
\\
& \le &
\frac{C}{(N_h+1)^2} \frac{1}{\ln(b/a)} \int_U |\nabla \tilde u|^2 \, dx.
\label{i0}
\end{eqnarray}
Let $v:\Omega_h\to\R^3$ be the function defined by
$$
v(x):= w(x_1, i_0hM+(-1)^{i_0} x_2, j_0hM+(-1)^{j_0} x_3) \qquad 
\hbox{for every } x\in\Omega_h.
$$
It is clear that $v\in W^{1,\infty}(\Omega_h;\R^3)$ and that it satisfies (\ref{trunc1}) by (\ref{wOm}). Moreover, since $\tilde u$ coincides with $u$ in $\Omega_h$ and it has been extended to $U$ by reflection, we have
\begin{equation} \label{i1}
\{x\in\Omega_h: \ u(x)\neq v(x)\} \subset \{ x\in S_{h,i_0j_0}:\ u(x)\neq w(x)\}
\end{equation}
and
\begin{equation} \label{i2}
\int_U |\nabla \tilde u|^2 \, dx \leq
(N_h + 2)^2 \int_{U_h} |\nabla \tilde u|^2 \, dx
\leq C(N_h + 2)^2 \int_{\Omega_h} |\nabla u|^2 \, dx,
\end{equation}
where the last inequality follows from \eqref{nablau} and \eqref{ud1}.
Now assertion (\ref{trunc2}) follows from (\ref{i0})--(\ref{i2}).

Finally, inequality \eqref{trunc3} is a standard consequence of \eqref{trunc2}.
This concludes the proof of Lemma~\ref{truncate}.
\end{proof}

\bigskip

\noindent {\bf Acknowledgments.}{ The authors gratefully acknowledge support from the Marie Curie
research training network MRTN-CT-2004-505226 (MULTIMAT). 
The first author was also partially supported by MIUR project 
``Calculus of Variations" 2004.}

\bigskip
\bigskip


\end{document}